\documentclass[final, 11pt]{article}
 
\pdfoutput=1

\usepackage[utf8]{inputenc}
\usepackage{mdwlist}
\usepackage{graphicx}
\usepackage{amsfonts}
\usepackage{ulem}
\usepackage{amsmath}
\usepackage{amsthm}
\usepackage{amssymb}
\usepackage{fancyhdr}
\usepackage{titlesec}
\usepackage{indentfirst}
\usepackage{booktabs}
\usepackage{verbatim}
\usepackage{color}
\usepackage{latexsym}
\usepackage{amscd}
\usepackage{esvect}
\usepackage{graphicx}
\usepackage{upgreek}
\usepackage{subcaption}   
\usepackage{caption}
\usepackage[page,header]{appendix}
\usepackage{titletoc}
\usepackage{mathrsfs}
\usepackage[numbers]{natbib}
\usepackage{float}
\usepackage{lipsum}
\usepackage[export]{adjustbox} 
\usepackage[linesnumbered,ruled,vlined, ruled,vlined]{algorithm2e}
\numberwithin{equation}{section}
\newtheorem{theorem}{Theorem}[section]

\newtheorem{remark}[theorem]{Remark}

\allowdisplaybreaks

\newcommand{\beq}{\begin{equation}}
\newcommand{\eeq}{\end{equation}}
\newcommand{\ben}{\begin{eqnarray}}
\newcommand{\een}{\end{eqnarray}}
\newcommand{\beno}{\begin{eqnarray*}}
\newcommand{\eeno}{\end{eqnarray*}}

\def\e{\varepsilon}

\begin{document}

\title{A deep neural network approach on solving the linear transport model under diffusive scaling}
\author{ Liu Liu\thanks{Department of Mathematics, The Chinese University of Hong Kong},~~
Tieyong Zeng \thanks{Department of Mathematics, The Chinese University of Hong Kong},~~
Zecheng Zhang\thanks{Department of Mathematics, Purdue University, West Lafayette, IN 47906, USA}}
\date{}
\maketitle

\begin{abstract}
    In this work, we propose a learning method for solving the linear transport equation under the diffusive scaling. 
    \textcolor{black}{Due to the multiscale nature of our model equation, the model is challenging to solve by using conventional methods}. We employ the physical informed neural network (PINN) framework, a mesh-free learning method that can numerically solve partial differential equations. 
    \textcolor{black}{Compared to conventional methods (such as finite difference or finite element type), our proposed learning method is able to obtain the solution at any given points in the chosen domain accurately and efficiently, which enables us to better understand the physics underlying the model.}
    In our framework, the solution is approximated by a neural network that satisfies both the governing equation and other constraints. 
    The network \textcolor{black}{is then trained with a combination of different loss terms.} 
    Using the approximation theory and energy estimates for kinetic models, we prove theoretically that the total loss vanishes as the neural network converges, upon which the neural network approximated solution converges pointwisely to the analytic solution of the linear transport model. 
    Numerical experiments for two benchmark examples are conducted to verify the effectiveness and accuracy of our proposed method. 
\end{abstract}

\section{Introduction}

The linear transport model, which describes the kinetic particles collision and absorption through a material medium while evolving in time, is an important physics model and arises in applications of various fields. 
In this paper, we study the linear transport equation under the diffusive scaling by using the deep neural network approach. 

The deep learning method has been successful in solving various problems in the computer science area; it improves the computation speed and the accuracy of many tasks such as the image classification, segmentation and language processing. Once the model is trained, one can use it by substituting the testing samples in the model, which usually is a forward process and only tensor productions are involved.
The data-driven solvers for partial differential equations (PDEs) have drawn an increasing attention due to their capability to encode the underlying physical laws governed by the model equation and give relatively accurate predictions for the unknowns. In particular, we mention the physical informed neural networks (PINNs) framework \cite{jin2020nsfnets, raissi2020hidden, raissi2019physics, kissas2020machine, sahli2020physics, han2018solving, wang2020understanding}. This method leverages the benefits of auto-differentiation of the current software and the underlying physics of the PDEs, which can be incorporated into the network by minimizing the losses composed of the PDE's residuals and other constraints such as the initial and boundary conditions. There are many other works developed, for example \cite{chung2020multirl, chung2020multi, zhang2020learning, Wang-Zhang20}. In this work, we apply this popular method to solve the important, hyperbolic-type kinetic equation. 

The Deep Learning method can resolve some difficulties in traditional finite difference-type numerical methods \cite{Karn-NN}, 
such as the expensive computational cost especially for problems with high-dimensional physical variables, the challenge dealing with complex boundary conditions, truncation of the velocity domain (in some models such as the Boltzmann equation, the velocity lies in the three-dimensional whole space thus a truncation is needed in numerical discretizations).  
Besides, the Deep Learning algorithm has the advantage of being intuitive and easy to be executed. For example, instead of designing mass (or momentum and energy, if applicable) conservative schemes for kinetic models, that are challenging in traditional numerical methods \cite{MM-Cons}, one can simply involve the derivative in time of conserved quantities of interests in the total loss function, as implemented in \cite{Hwang}. 

We mention some other advantages of using the DNN approach to solve the linear transport model: 
a) to obtain the distribution function at any given $t$, $x$, $v$, instead of getting only discrete values at the uniform mesh in traditional finite volume or finite element numerical methods; as a mesh-free method, works efficiently for high-dimensional physical space problems; 
b) to avoid high computational cost on the simulation of kinetic equations due to the velocity variable $v$, and the integral-based, nonlocal collision operators appeared in other complicated kinetic models such as the Boltzmann or Landau equations. 
Nevertheless, we mention that there are indeed some weaknesses of the Deep Learning approach.
First, there is no guarantee that the Deep Learning algorithm will converge and it is practically difficult to show their convergence. It is also hard to evaluate the accuracy of the DNN approach in contrast with traditional numerical methods.

This paper is organized as the following. In Section \ref{sec:LTE}, we introduce the background of the linear transport model under the diffusive scaling. 
In Section \ref{sec:NN}, we review and discuss the neural network framework and method for solving general PDEs. 
Two main convergence results are given in Section \ref{sec:Anal}, showing that 
1) the loss function goes to zero as the neural network converges; 
2) the neural network solution converges point-wisely to the analytic solution when the loss function converges to zero. 
The effectiveness and accuracy of our proposed method, including the choice of weights in the total loss function will be presented in Section \ref{sec:Num}. 
Finally, we summarize the paper and mention some future work in Section \ref{sec:Con}.

%--------------------------------------
\section{The linear transport model}
\label{sec:LTE}

The linear transport model, which describes the kinetic particles collision and absorption through a material medium while evolving in time, arises in many applications, such as atmosphere and ocean modeling \cite{coakley2014atmospheric,stamnes2017radiative,zdunkowski2007radiation}, astrophysics \cite{peraiah2002introduction} or nuclear physics. Such problems usually involve several orders of 
magnitudes of length scales characterized by the Knudsen number, defined as the ratio of the mean free path over a typical length scale such as the size of the spatial domain. 

We consider the linear transport equation under diffusive scaling, with one-dimensional space and velocity variable ($x\in\Omega\subset\mathbb R$, $v\in [-1,1]$) and given by
\begin{align}
\label{LTE}
\begin{split} &\e\partial_t f + v\cdot\nabla_x f = \frac{1}{\e} \mathcal L(f), \qquad f(t=0,x,v)=f_0(x,v), \\[4pt]
& \mathcal L(f) = \sigma(x)\left[\frac{1}{2}\int_{-1}^1 f(t,x,v)dv - f(t,x,v)\right].
 \end{split}
 \end{align}

%--------------------------------------
\section{The Neural Network approach}
\label{sec:NN}
We briefly review the Deep Neural Network (DNN) structure and approach introduced in \cite{Hwang}, where the 
one-dimensional kinetic Fokker-Planck equation is studied. 
Denote the approximated function as $f^{nn}(t,x,v,z;m,w,b)$ and suppose DNN has $L$ layers; the input layer takes $(t,x,v)$ as input 
and the final layer gives $f^{nn}(t,x,v,z;m,w,b)$ as the output. The relation between the $l$-th and $(l+1)$-th layer ($l=1,2, \cdots L-1)$ is 
given by
$$ u_{j}^{(l+1)}=\sum_{i=1}^{m_{l}} w_{j i}^{(l+1)} \bar{\sigma}_{l} (z_{i}^{l})+b_{j}^{(l+1)}, $$
where $m=\left(m_{0}, m_{1}, m_{2}, \dots, m_{L-1}\right)$, the weights $w=\left\{w_{j i}^{(k)}\right\}_{i, j, k=1}^{m_{k-1}, m_{k}, L}$ and the bias $b=\left\{b_{j}^{(k)}\right\}_{j=1, k=1}^{m_{k}, L}$
are given in \cite{Hwang}, which we refer to for details. 

Regarding the optimization algorithm, we use {\it Adam} optimization algorithm, which is an extended algorithm of the stochastic gradient descent and is popularly used in the applications of the deep learning.

\subsection{Definition of loss functions}

The loss function for the governing linear transport equation \eqref{LTE} is defined by
\begin{align}
\label{Loss1}
\begin{split}
&\displaystyle\quad Loss_{GE} \\[4pt]
&\displaystyle  = \int_0^T \int_{\Omega}\int_V \Big| \e\partial_t f^{nn}(t,x,v;m,w,b) + v_k \partial_x f^{nn}(t,x,v;m,w,b)   \\[4pt]
&\displaystyle\quad - \frac{1}{\e} \mathcal L(f^{nn})(t,x,v;m,w,b) \Big|^2\,  dv dx dt \\[4pt]
&\displaystyle \approx \frac{1}{N_{i,j,k}} \sum_{i,j,k} \Big| \partial_t f^{nn}(t_i, x_j, v_k; m,w,b) + v_k \partial_x f^{nn}(t_i, x_j, v_k; m,w,b) \\[4pt]
&\displaystyle \quad - \frac{1}{\e} \mathcal L(f^{nn})(t_i, x_j, v_k; m,w,b)\Big|^2, 
\end{split}
\end{align}
where $N_{i,j,k} = N_i N_j N_k$, and the collision operator
$$\mathcal L(f^{nn})(t_i, x_j, v_k; m,w,b) = \sigma(x_j)\left[\frac{1}{2 N_k} \sum_k f^{nn}(t_i, x_j, v_k; m, w, b) - f^{nn}(t_i, x_j, v_k; m,w,b)\right]. $$

Define $n_x$ the {\it unit} outward normal vector on the boundary $\partial\Omega$, and $\gamma \stackrel{\text { def }}=\partial\Omega\times [-1,1]$. 
This phase boundary can be split into an outgoing boundary $\gamma_{+}$, incoming boundary $\gamma_{-}$ and a singular boundary $\gamma_0$, defined by
\begin{align*}
\gamma_{+} \stackrel{\text { def }}=\left\{(x, v) \in \partial\Omega \times [-1,1]: n_x \cdot v>0\right\}, \\[2pt]
\gamma_{-} \stackrel{\text { def }}=\left\{(x, v) \in \partial\Omega \times [-1,1]: n_x \cdot v<0\right\}, \\[2pt]
\gamma_{0} \stackrel{\text { def }}=\left\{(x, v) \in \partial\Omega \times [-1,1]: n_x \cdot v=0\right\}. 
\end{align*}

We now define the loss terms for the initial condition and the boundary conditions: 
\begin{align}
\label{Loss2}
\begin{split}
\displaystyle Loss_{IC} &= \int_{\Omega}\int_V \left|f^{nn}(0,x,v) - f_0(x,v)\right|^2\, dv dx\\[4pt]
\displaystyle & \approx \frac{1}{N_{j,k}}\sum_{j,k}\left| f^{nn}(0,x_j,v_k) - f_0(x_j,v_k)\right|^2. 
\end{split}
\end{align}
For inflow boundary condition $f(t,x,v)|_{\gamma^{-}}=g(t,x,v)$ on $x\in\partial\Omega$, %and for each fixed $z$ (similarly for other types), 
one has 
\begin{align}
\label{Loss3}
\begin{split}
\displaystyle Loss_{BC} &= \sum_{x\in \partial\Omega} \int_0^T \int_V  \left|f^{nn}(t,x,v;m,w,b) - g(t,x,v;m,w,b)\right|^2\, dv dt  \\[4pt]
\displaystyle &\approx \frac{1}{|\partial\Omega| N_{i,k}} \sum_{i,k} \left| f^{nn}(t_i,x,v_k;m,w,b) - g(t_i, x, v_k;m,w,b)\right|^2, 
\end{split}
\end{align}
where $|\partial\Omega|$ denotes the volume of the spatial boundary. 
Add up \eqref{Loss1}, \eqref{Loss2} and \eqref{Loss3}, with appropriate weights 
$\{\lambda_g, \lambda_i, \lambda_b\}$, the total loss function is defined by
\begin{equation} Loss_{Total} = \lambda_g\, Loss_{GE} + \lambda_i \, Loss_{IC} + \lambda_b \, Loss_{BC}. 
\end{equation}

%-------------------------------
\section{Analysis results}
\label{sec:Anal}

In this section, we show the two main theoretical results. 

Recall that neural network architecture was first introduced in \cite{NN1943}. Later in \cite{UAT}, Cybenko established sufficient conditions where a continuous function can be approximated by finite linear combinations of single hidden layer neural networks, followed by the work in \cite{Multilayer} that extends the theory to the multi-layer network case. 
We paraphrase the Universal Approximation Theorem \cite{UAT} to the form that is needed in our context: 
If $f$ solves the linear transport model \eqref{LTE} and is sufficiently smooth in all variables, then there exists a two-layer neural network
$$f^{nn}(t, x, v)=\sum_{i=1}^{m_1} w_{1 i}^{(2)} \bar{\sigma}\left(\left(w_{i 1}^{(1)}, w_{i 2}^{(1)}, w_{i 3}^{(1)}\right) \cdot(t, x, v)+b_{i}^{(1)}\right)+b_{1}^{(2)}$$
such that 
$$ || f - f^{nn} ||_{L^{\infty}(K)} < \eta, \quad
|| \partial_t (f - f^{nn}) ||_{L^{\infty}(K)} < \eta, \quad
|| \nabla_x (f - f^{nn}) ||_{L^{\infty}(K)} < \eta,
$$
where the domain $K = [0, T] \times\Omega \times [-1,1]$. 
%$\left\|D^{k} f-D^{k} f^{nn}\right\|_{L^{\infty}([0, T] \times\Omega \times [-1,1])}<\eta$, where $D^k$ denotes the $k$-th order derivatives. 
Note that our analysis results in this section apply to the neural network with multiple layers, due to \cite{Multilayer}. For simplicity, we review the Universal Approximation Theorem in its two-layer form. 

\subsection{Convergence of the loss function} 

We first show that a sequence of neural network solutions to \eqref{LTE}
exists such that the total loss function converges to zero.  
\begin{theorem}
\label{Thm1}
Let $f$ solves the equation \eqref{LTE} and satisfies $f \in C^1([0,T])\cap C^1(\Omega\times V)$. 
Then there exists a sequence of neural network parameters $\{m_{[j]}, w_{[j]}, b_{[j]}\}_{j=1}^{\infty}$ such that the sequence of DNN solutions with $m_{[j]}$ nodes, given by $\{f_j(t,x,v) = f^{nn}(t,x,v;m_{[j]}, w_{[j]}, b_{[j]})\}_{j=1}^{\infty}$ satisfies 
$$ Loss_{\text{Total}}(f_j)\to 0, \qquad \text{as  } j\to\infty. $$
\end{theorem}

{\textbf{Proof.} }
Denote the velocity domain $\mathcal D = [-1,1]$ and $\mu(V)=\int_{\mathcal D} dv = 2$. 
Define \begin{equation}\label{d-GE} d_{ge, j}(t,x,v,z) := - \left[\e \partial_t + v\cdot\nabla_x\right] f_j + \frac{1}{\e} \mathcal L(f_j). \end{equation}
Integrate $ | d_{ge,j} |^2$ over $[0,T]\times\Omega\times\mathcal D$, which is equivalent to
\begin{equation}\label{Eqn} |\e \partial_t ( f - f_j) + v\cdot \nabla_x ( f - f_j) - \frac{1}{\e}\mathcal L(f) + \frac{1}{\e} \mathcal L(f_j) |^2 \end{equation}

We will show that the loss term \eqref{Loss1} is bounded by $O(\eta^2)$. 
For the first two terms in \eqref{Eqn}, it is bounded by $\eta^2\, T \mu(V)\, |\Omega|\, (\e^2 + \mu(V)^2)$. 
Observe that for any function $g$,  \begin{equation}\label{L-g} \int_{-1}^1 \left({\mathcal L}(g)\right)^2 dv =  -\frac{\sigma^2}{2}\left(\int_{-1}^1 g(v)dv \right)^2 +  \sigma^2 \int_{-1}^1 g^2(v)dv \leq \sigma^2\int_{-1}^1 g^2(v)dv.   \end{equation}
Since \begin{equation}\label{f-L} ||f-f_j||_{L^2(\mathcal D)}^2 \leq ||f-f_j||_{L^{\infty}(\mathcal D)}^2\, \mu(V) < C \eta^2, \end{equation}
Let $g=f-f_j$ in \eqref{L-g} and by \eqref{f-L}, we have 
 $\| \mathcal L(f) - \mathcal L(f_j) \|_{L^2_v}^2 = \| \mathcal L(f-f_j) \|_{L^2_v}^2 < O(\eta^2)$. 
So the last two terms in \eqref{Eqn}: 
 \begin{equation}\label{L-F} \left\| \frac{1}{\e}\mathcal L(f) - \frac{1}{\e}\mathcal L(f_j) \right \|_{L^2(\Omega\times\mathcal D)}^2 <  
O \left(\frac{\eta^2}{\e^2}\right), 
\end{equation}
since $\Omega$ is bounded. 
Therefore, the lost term $Loss_{GE}$ in \eqref{Loss1} bounded by $O\left(\e^2 \eta^2 + (\frac{\eta}{\e})^2\right)$.

\hspace{10cm}

For the inflow boundary condition, $Loss_{BC}$ is bounded by 
$$ || f_j - f ||_{L^2(\gamma_{T}^{-})}^2 \leq T \mu(V) |\partial\Omega| \, || f_j - f ||_{L^{\infty}(\gamma_{T}^{-})}^2 
\leq  T \mu(V) |\partial\Omega|\, || f_j - f ||_{L^{\infty}([0,T]\times\Omega\times \mathcal D)}^2 \leq O(\eta^2), $$
where $\gamma_{T}^{\pm} := [0,T]  \times \gamma_{\pm}$. Note that the specular reflection boundary condition works similarly. 
For the initial data, denoted by $f_{j, in}$ and $f_{in}$ respectively, 
$$ Loss_{IC} = || f_{j, in} - f_{in} ||_{L^2(\Omega\times \mathcal D)}^2 \leq || f_{j, in} - f_{in} ||_{L^{\infty}(\Omega\times \mathcal D)}^2\,  |\Omega|\, \mu(V)
\leq O(\eta^2). $$
Set $\eta= \eta_j = \frac{1}{j}$, combine all the loss terms \eqref{Loss1}-\eqref{Loss3}, we conclude that
$$ Loss_{Total}(f_j) \leq O(\frac{1}{\e^2 j^2}). $$
Therefore, $Loss_{Total}(f_j)\to 0$, as $j\to\infty$.

%-------------------------------------------------------------------------------------------

\subsection{Convergence of the neural network approximated solution} 
In this section, we show that with the parameters $\{m_{[j]}, w_{[j]}, b_{[j]}\}_{j=1}^{\infty}$ equipped, the neural network in Theorem \ref{Thm1} converges to the analytic solution of the linear transport model. 
\begin{theorem}
Let $\{m_{[j]}, w_{[j]}, b_{[j]}\}_{j=1}^{\infty}$ be a sequence defined in Theorem \ref{Thm1} and $f$ solves the linear transport model \eqref{LTE}. As a consequence, $Loss_{\text{Total}}(f_j)\to 0$ implies that
$$ ||f_j(\cdot, \cdot, \cdot, m_{[j]}, w_{[j]}, b_{[j]}) - f||_{L^{\infty}([0,T]; L^2(\Omega\times\mathcal D))}
\to 0, $$
for some finite time $t\in[0,T]$, physical variables $(x, v)\in \Omega\times \mathcal D$. 
\end{theorem}

{\textbf{Proof. }}
Recall the definition \eqref{d-GE}, then 
\begin{equation}\label{GE} \left[\partial_t + v\cdot\nabla_x \right] \{f-f_j\} = d_{ge,j}(t,x,v) + \frac{1}{\e}\mathcal L(f) - \frac{1}{\e}\mathcal L(f_j). \end{equation}
Define $$d_{ic,j}(x,v) := f_0(x,v) - f_j(0,x,v), $$ in addition to 
 $$ d_{bc,j}(t,x,v) := g(t,x,v)- f_j(t,x,v)\quad \text{  at     }\, (t,x,v)\in\gamma^{-}_{T}, $$ for inflow boundary condition, and 
$d_{bc,j}(t,x,v) :=f_j(t,x,-v)-f_j(t,x,v)$ for specular reflection boundary condition. 

The $L^2$ norms and and inner products below stand for $|| \cdot ||_{L^2(\Omega\times\mathcal D)}$. 
Multiplying $(f-f_j)$ onto \eqref{GE} and integrating over $\Omega\times \mathcal D$, one gets
\begin{align}
\label{GE1}
\begin{split}
&\displaystyle\quad \int_{\Omega}\int_{\mathcal D} \e\partial_t (f-f_j)^2\, dvdx + \int_{\gamma^{+}}(f-f_j)^2\, v\cdot n_x\, d\gamma - 
\int_{\gamma^{-}}d_{bc,j}^2\, |v\cdot n_x| \, d\gamma \\[4pt]
&\displaystyle = \frac{2}{\e} \big\langle \mathcal L(f - f_j), f-f_j \big\rangle_{L^2} + 2 \big\langle d_{ge,j}, f-f_j \big\rangle_{L^2}. 
\end{split}
\end{align}
Recall \eqref{L-F} that $\|\mathcal L(f - f_j)\|_{L^2}^2 \leq C \eta^2$. 
The right-hand-side in \eqref{GE1} is bounded by
$$ \frac{2}{\e} \big\langle \mathcal L(f - f_j), f-f_j \big\rangle_{L^2} \leq \frac{1}{\e} \| \mathcal L(f - f_j)\|_{L^2}^2 + \frac{1}{\e}\| f - f_j \|_{L^2}^2, $$
and since $\int_{\gamma^{+}}(f-f_j)^2\, v\cdot n_x\, d\gamma \geq 0$, thus
$$ \e \frac{d}{dt}||f-f_j||_{L^2}^2 \leq \underbrace{\int_{\gamma^{-}}d_{bc,j}^2\, |v\cdot n_x| \, d\gamma + ||d_{ge,j}||_{L^2}^2 + \frac{ C \eta^2}{\e}}_{:=H(t)} + (1+\frac{1}{\e}) ||f-f_j||_{L^2}^2. $$
Since $\int_{\gamma^{-}}d_{bc,j}^2\, |v\cdot n_x| \, d\gamma \leq \int_{\gamma^{-}}d_{bc,j}^2\, d\gamma$, and the definitions that
$$ \int_0^t  \int_{\gamma^{-}}\left(d_{bc,j}(s,\cdot, \cdot, \cdot)\right)^2\, d\gamma ds = Loss_{BC}, \quad 
\int_0^t  ||d_{ge,j}(s, \cdot, \cdot, \cdot)||_{L^2}^2\, ds = Loss_{GE}, $$ then 
\begin{align*}
\notag\displaystyle ||f-f_j||_{L^2(\Omega\times \mathcal D)}^2 &\leq e^{(\frac{1}{\e} + \frac{1}{\e^2})t} Loss_{IC} +  \frac{1}{\e} e^{(\frac{1}{\e} + \frac{1}{\e^2})t} \int_0^t H(s) ds \\[4pt]
&\notag\displaystyle \lesssim e^{(\frac{1}{\e} + \frac{1}{\e^2})t}\left(Loss_{IC} + \frac{1}{\e} Loss_{BC} + \frac{1}{\e} Loss_{GE} +  C_0\, \frac{\eta^2}{\e^2} t \right) \\[4pt]
&\displaystyle \leq e^{(\frac{1}{\e} + \frac{1}{\e^2})t} \left(\frac{1}{\e}Loss_{Total} + C_0\, \frac{\eta^2}{\e^2} t \right). 
\end{align*}

We already know from Theorem \eqref{Thm1} that $Loss_{Total}(f_j) \leq O(\frac{\eta^2}{\e^2})$,  
now take $L^{\infty}$ in $t\in [0,T]$ with $\eta=\eta_j=\frac{1}{j} \to 0$, then 
\begin{equation}\label{R1} ||f-f_j||_{L^2(\Omega\times\mathcal D)}^2 \leq C^{\prime}  e^{(\frac{1}{\e} + \frac{1}{\e^2})t} \,\eta^2 \left(\frac{1}{\e^3} + C_0 \frac{t}{\e^2}\right). 
\end{equation}
Therefore,  $$||f-f_j||_{L^{\infty}([0,T]; L^2(\Omega\times\mathcal D))} \to 0, \quad \text{as   }\, j\to \infty\,. $$

\begin{remark}
In \cite{Liu-Jin18}, Liu et al. employed the hypercoercivity and Lyapunov-type functionals to conduct sensitivity analysis for a general class of kinetic equations with random uncertainties and multiple scales. It is possible to adopt that framework and energy estimates analysis to study the convergence of the DNN solution (or even improve to a stronger result of exponential decaying in time), as presented in this section. However, due to its high complexity and main focus of the current manuscript, we defer to a future work. 
Also, we mention that a similar analysis techniques of this section can be found in \cite{CLM21}, where the radiative transfer model in the kinetic regime ($\e=1$), with numerical examples of different applications and its long time behavior was studied by the DNN approach. 
\end{remark}
%---------------------------------------------------------------------------------------------------------------
\section{Numerical implementation}
\label{sec:Num}
In this section, we are going to verify the proposed method with two examples. 
\textcolor{black}{The first example aims to show the accuracy of the method, and a more challenging problem with a boundary layer is studied in the second test, which verifies that our proposed method is able to capture the physical properties of the model.}
We will compute the solutions of the model problems by using the numerical method shown in the Appendix, then compare it with the proposed deep learning method. As discussed before, the deep learning method is a mesh-free method, hence inherits all the benefits of the mesh free numerical methods meanwhile keeps a low computational cost in the testing phase. We will first introduce the grid points of the testing data which are used in the numerical computational. Details of the two examples will then be presented.

\subsection{Reference solution grid points and set-up}
\label{sec::ref}
\textcolor{black}{In this section, we are going to show the computation details of the reference solution.}
In order to obtain the reference solutions from the conventional solver, we adopt the robust, implicit yet explicitly implementable Asymptotic-Preserving (AP) numerical scheme in the even-odd decomposition framework \cite{JPT2}, which works with high resolution uniformly with respect to the scaling parameter $\varepsilon$. We will use solutions of this solver as the reference solutions to compare with the DNN approximated solutions in our numerical tests. For the convenience of readers, we review this scheme in the Appendix. 
 
Let $\Omega = [0,1]$. The grid points of $t$ and $x$ for the training are chosen uniformly as follows: 
$$ \{ (t_i, x_j) \}_{i,j} \in [0,T]\times [0,1], \quad\text{  with fixed  } \Delta t,\, \Delta x.  $$
 The integral in velocity can be computed using the quadrature rule, only few points are needed, such as $N_v=32$ with $\{ v_k \}_{k=1}^{N_v}$. 
Use the grids $\{ (t=0, x_j, v_k) \}_{j,k}$ for the initial condition and 
$\{ (t_i, x=0 \text{ or } 1, v_k)\}_{i,k}$ for the boundary condition. 
We save the data of $f$ at fixed $t$ and $x$, for all velocity points at $[-1,1]$, chosen uniformly with $N_v=32$. 
In the conventional finite difference AP solver, set $x_i$ ($1\leq i \leq N$) and $N=40$. Here $\Delta x = 0.025$. 
%The $N$ interior spatial points (besides the ghost points outside the boundary) are: \\ 
%$[\frac{\Delta x}{2}: \Delta x: 1-\frac{\Delta x}{2}]$. 
According to the parabolic CFL of the AP scheme, $\Delta t \sim O((\Delta x)^2)$. In our test, $\Delta t= 0.5 (\Delta x)^2$. %$\Delta t=3.1250e-04$. 
 \textcolor{black}{To summarize, in both of our experiments, we test our models on the reference solutions which are evaluated basing on the following fine mesh:}
\begin{equation*}
\{ (t_i, x_j, v_k)\}_{i,j,k} \in [0,T]\times [0,1]\times[-1,1], \qquad \Delta t=0.5\Delta x^2,\, \Delta x=\frac{1}{40}, \,\Delta v=2/(33-1). 
\end{equation*}

\subsection{Test I}
\label{sec::test2}
%\textcolor{red}{Liu: can you introduce the background of the problem and why we need to test on this example?}
We consider a benchmark test for studying the linear transport model \cite{JPT2}.  %\textcolor{red}{Liu, can you introduce a little bit about the difficulties in solving this problem numerically.} \textcolor{cyan}{This is actually the simplest test, without much difficulties given these initial and boundary conditions. However, I add a remark for general challenges (we do not test those examples though, i.e., when the Knudsen number is not a constant and varies dramatically with space. So I would not want to emphasize too much here.)}
\textcolor{black}{In practice, the scaling parameter (mean free path) may differ in several orders of magnitude from the rarefied regimes to the diffusive regimes in one problem, thus developing methods that work {\it uniformly} with respect to this parameter becomes important. Our learning method manages to achieve this goal, one can test on models at any given Knudsen number without bringing additional challenges. In this test, we consider a smooth initial condition and the diffusive scaling.}
Let the initial distribution be the double-peak Maxwellian: 
\begin{align} 
\label{IC}
\left\{
\begin{array}{l}
\displaystyle \rho_{0}(x) = 1 + \frac{1}{2}\sin(2\pi x), \\[10pt]
\displaystyle T_{0}(x) = \frac{5+2\cos(2\pi x)}{20}, \\[10pt]
 \displaystyle f(t=0,x,v)=\rho_0\left[\exp\left(-\left(\frac{v-0.75}{T_0}\right)^2\right)+\exp\left(-\left(\frac{v+0.75}{T_0}\right)^2\right)\right].
  \end{array}\right.
\end{align}
Periodic boundary condition is considered. Assume the scattering coefficient $\sigma=1$ and $\varepsilon=10^{-2}$. 
%{\color{blue} Filename meaning: EX1 or EX2 means Test I or Test II in the numerical tests section. \\
%The index followed `Time' means number of time steps. For example, `EX1Time8' means Test I, the distribution $f$ as a matrix of $x$ and $v$. 
%For the `.txt' file, it means that the first $32$ numbers are $f$ values at the first point of the spatial mesh, i.e., $\frac{\Delta x}{2}$ with $32$ values of velocities; 
%the second $32$ numbers are $f$ values at the 2nd point of the spatial mesh, i.e., $\frac{3\Delta x}{2}$ with all $32$ values of velocities, so on so forth. 
%The two files `Velocity.txt' and `Weights.txt' are given as two separate files. }

Similar to \cite{Hwang}, we make the data of grid points for each variable for the DNN training. We generate the mesh uniformly with
$\Delta t=0.5\Delta x^2,\, \Delta x=1/20, \,\Delta v=2/(33-1)$. 
This grid points setting provides us with $17000$ samples associated with the the governing equation loss (\ref{Loss1}), $850$ samples associated with the boundary loss (\ref{Loss3}) and $340$ samples associated with the initial loss (\ref{Loss2}). We train the model which will be detailed later and test the performance on a pre-generated data set which is specified in the section (\ref{sec::ref}).

\begin{figure}[H]
\centering
\mbox{
\includegraphics[scale = 0.3]{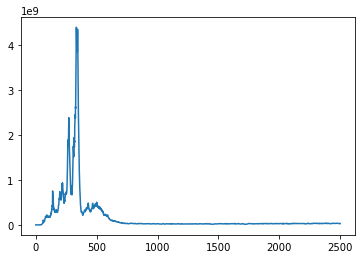}
\includegraphics[scale = 0.3]{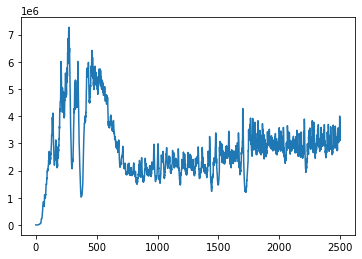}
}
\caption{y axis: weight; x axis: training epoch. Left: boundary loss; right: initial condition loss. The x- axis is the training epochs and y-axis is the weight.}
\label{weight_exp2}
\end{figure}

We use a 4-layers full connected network ($3\times 256\rightarrow 256\times 256\rightarrow 256\times 256\rightarrow 256\times 1$) with bias and activated by Tanh function. the network is trained with the Adam gradient descent with 2500 epochs. The initial learning rate is set to be 0.005 and we used a step scheduler with step size 750 and the reduce rate 0.95.
\textcolor{black}{
The PINNs model is not easy in the model training \cite{wang2020understanding, wang2020and} and there is no reliable understanding or general method to tackle this issue \cite{wang2020eigenvector}. One of the difficulties is the weight tuning.
\textcolor{black}{In particular, when $\epsilon$ is small, the governing equation loss (\ref{Loss1}) will dominate the total loss and there are multiple scales in the loss function, bringing us an additional difficulty in the weight tuning.}
Regarding this, we are going to adapt the strategy in \cite[Algorithm 1]{wang2020understanding}}. To make the current work more readable, we cite and present their algorithm in the Appendix. 
The weights for initial and boundary losses are depicted in Figure (\ref{weight_exp2}). Note that the weight for the general equation is normalized to $1$. \textcolor{black}{Our large scale experiments show that using the above-mentioned adaptive weights gives us more accurate prediction results than the pre-set constant weights.}

\begin{figure}[H]
\centering
\mbox{
\includegraphics[scale = 0.4]{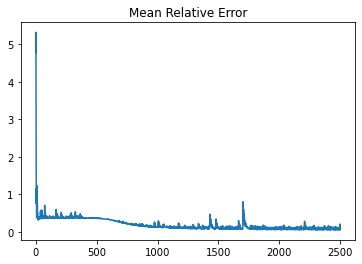}
}
\caption{Relative loss defined by \eqref{relative_loss}, with respect to the number of training epochs. 
y axis: relative loss as defined in (\ref{relative_loss}); x axis: training epoch. }
\label{exp2_loss}
\end{figure}

\begin{figure}[H]
\centering
\mbox{
\includegraphics[scale = 0.3]{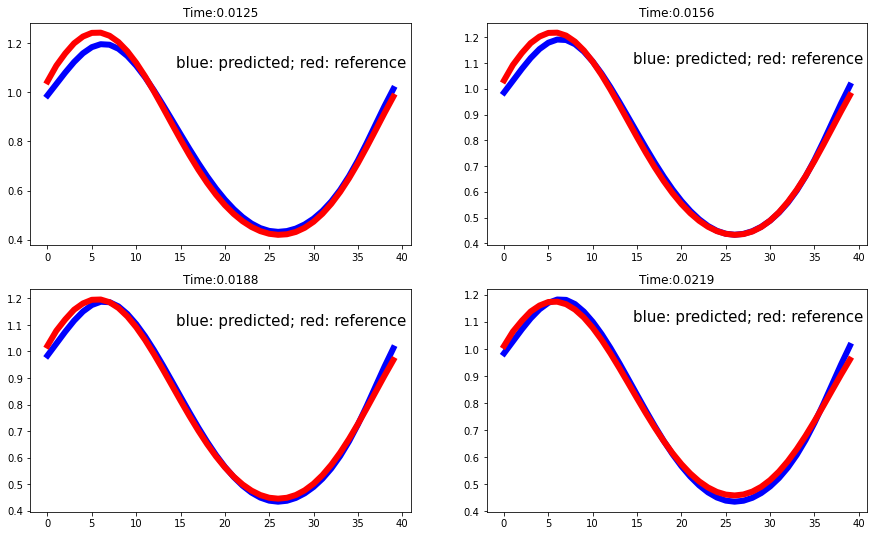}
}
\caption{Plot of the density at different output time. 
The blue curve is the NN predicted solution and the red curve is the reference solution. }
\label{exp2_results}
\end{figure}

We computed the relative error of (\ref{r_integral}), more precisely:
$$
\text{relative error} = \frac{\|u_{\text{True}}-u_{\text{Predicted}} \|}{\|u_{\text{True}}\|},
\label{relative_loss}
$$
where $u_{\text{True}} = \int_{-1}^{+1}f dv$ and we denote the true solution as $f$. The evolution of the relative errors with respected to the training epoch is shown in Figure (\ref{exp2_loss}).
\textcolor{black}{Note that the largest relative error stabilizes at less than 10 percents when the neural network converges and the number of epochs reach about $2000$.} %\textcolor{blue}{Liu, do we need to mention this since 10 is not small error? Also, do we need to include the relative error figure?}
%\textcolor{cyan}{I think it is OK now, we can just leave like this and see what the reviewers get to say...}
We also plot the (\ref{r_integral}) of the true and predicted solution in Figure (\ref{exp2_results}) at different time.
\textcolor{black}{We can observe from the Figure (\ref{exp2_results}) that the proposed neural network method gives an accurate approximation to the reference solution; this approximation is close even when time is large. Also, periodic boundary conditions are satisfied. This example verifies that our method is accurate and then we can benefit from the many advantages of the learning algorithms. One of the advantages is the high efficiency in predicting new samples once the model is trained. This is very useful since our framework is mesh-less and compare to the conventional numerical methods, we are able to calculate the solution at any given points in the support of the equation. One hence can study the physics of the model efficiently. }

\subsection{Test II}
\begin{figure}[H]
\centering
\mbox{
\includegraphics[scale = 0.3]{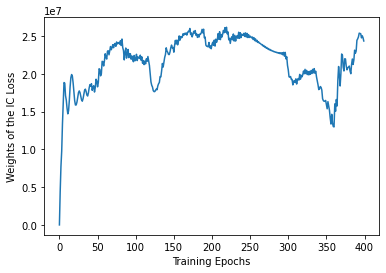}
\includegraphics[scale = 0.3]{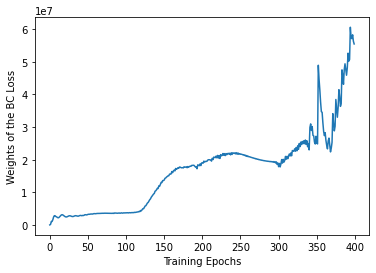}
}
\caption{y axis: weight; x axis: training epoch. Right: boundary loss; 
left: initial condition loss. The loss for the governing equation is normalized to 1.}
\label{weight_exp1}
\end{figure}

In this experiment, we consider another benchmark test with the incoming boundary data \cite{JPT2}. This problem is more complicated and challenging, since the solution contains the boundary layer, and we manage to see the DNN approach can capture the solution behavior especially near the boundary. 
%\textcolor{red}{Liu: can you introduce the background of the problem and why we need to test on this example? Why this example is hard to do?}
The initial condition is given by
$f(0, x, v)=0$, and boundary conditions are
$$f(t, 0, v)=1, \quad v \geq 0;  \quad f(t, 1, v)=0, \quad v \leq 0. $$
Consider the diffusive regime with $\varepsilon=10^{-3}$, and $\sigma=1$. Set $\Delta t = 0.5\cdot\Delta x^2$ and $\Delta x = 1/25$. 
This gives us 41225 samples for the governing equation loss, 1746 samples for the boundary loss and 425 for the initial condition loss. We use a 4-layers fully connected network ($3\times 256\rightarrow 256\times 256 \rightarrow 256\times 256 \rightarrow 256\times 1$) activated by Tanh and trained by the Adam gradient descent for 400 epochs. The initial learning rate is set to be 0.0005 and a step schedule is used which will reduced the rate by 5 percents every 50 epochs. Same as the first example, we will adopt the adaptive weight balance strategy. The weight evolution of the initial and boundary loss is shown in Figure (\ref{weight_exp1}).

\begin{figure}[H]
\centering
\mbox{
\includegraphics[scale = 0.3]{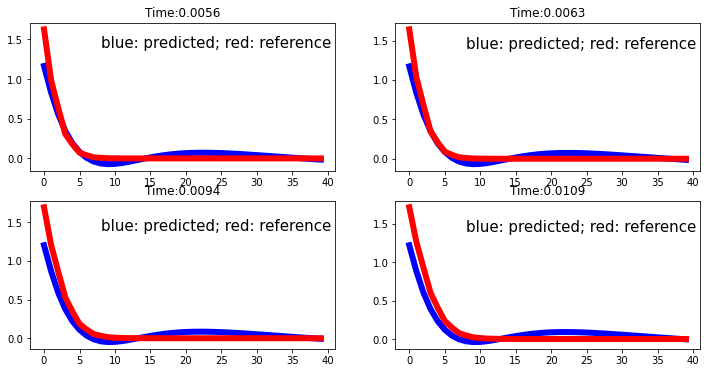}
}
\caption{y axis: can be calculated by (\ref{r_integral}); x axis: points in spatial direction x. Each figure is the result at a time step, please see the subtitles for the details of the time level. The blue curve is the predicted solution and the red curve is the reference solution. }
\label{exp1_results}
\end{figure}

We finally present the results of the experiment in figure \ref{exp1_results}.
\textcolor{black}{The most challenging property of this example is the existing of the boundary layer. We can see from the Figure (\ref{exp1_results}), the proposed method is able to capture the boundary layer. This is the consequence of minimizing the combined loss and the solution hence has the all the physical properties of the equation.}

%----------------------------------
\section{Conclusion and future work}
\label{sec:Con}

In this work, we proposed to solve the linear transport model by a learning method. 
%{\color{red} In real world applications, the diffusion scaling is usually small and brings the multiscale property to the model. 
%The conventional method is hard to capture the multiscale features and we hence propose the learning method. }
%{\color{cyan}The above sentences are not technically right, there is no such thing for our model...I rewrite it as the following. }
{\color{black} We consider the diffusive scaling with a small Knudsen number in our numerical tests, while our analysis applies to all orders of the Knudsen number. The asymptotic-preserving solver \cite{JPT2}, as a robust traditional numerical method, was designed to tackle the stiffness of the model brought by the small relaxing parameter, without resolving the numerical discretizations. However, it requires a solid understanding of kinetic theory and not very feasible to practical applications in physics or engineering. 
Our learning method, on the other hand, is mesh-free, easy to implement (no matter how complicated the initial or boundary conditions are), provides the numerical solution at any given points, while does not need a strong background on kinetic theory, thus making it more applicable to general research fields.} Theoretically, we prove that the total loss function vanishes as the network converges, then show that the neural network solution converges to the analytic solution pointwisely. 
In order to demonstrate the advantages of the learning method, we test on two benchmark examples, whose results show that our method is accurate and can capture the quantities of interests accurately, given challenging initial or boundary conditions. 

In the future, we will extend our proposed method to high-dimensional kinetic problems with uncertainties, and develop new training methods--in particular work on weight balancing of different loss terms. We may also consider applying the PINN framework to solve inverse problems associated with kinetic models. 

\bigskip
\bigskip

%----------------------------------
{\bf{\Large Appendix}}

\hspace{2cm}

{\bf A. The asymptotic-preserving method }
We briefly recall \cite{JPT2} for the reformulation to diffusive relaxation system of the linear transport equation \eqref{LTE}, and its diffusion limit system 
as $\varepsilon\to 0$. This also prepares us to study the asymptotic behavior of the distribution function, which will be studied in a follow-up work. 

First, we split \eqref{LTE} into two equations for $v>0$: 
\begin{equation}
\label{split_LTE}
\begin{split}
\displaystyle 
&\varepsilon \partial_t  f(v) + v \partial_x f(v) = \frac{\sigma(x)}{\varepsilon}\left(\frac{1}{2}\int_{-1}^1f(v)\, dv-f(v)\right), \\[6pt]
\displaystyle 
&\varepsilon \partial_t f(-v) - v \partial_x f(-v) = \frac{\sigma(x)}{\varepsilon}\left(\frac{1}{2}\int_{-1}^1f(v)\, dv-f(-v)\right), 
\end{split}
\end{equation}
In this case consider the even and odd parities 
\begin{equation*}
\begin{split}
\displaystyle  r(t,x,v) &= \frac{1}{2}[f(t,x,v) + f(t,x,-v)], \\[4pt]
\displaystyle  j(t,x,v) &= \frac{1}{2\varepsilon}[f(t,x,v) - f(t,x,-v)].
\end{split}
\end{equation*}
Adding and subtracting the two equations in \eqref{split_LTE} leads to 
\begin{equation}
\label{RJ}
\left\{
\begin{split}
\displaystyle  &\partial_t r + v \partial_x j = \frac{\sigma(x)}{\varepsilon^2}(\rho-r),  \\[4pt]
\displaystyle  &\partial_t j + \frac{v}{\varepsilon^2} \partial_x r = -\frac{\sigma(x)}{\epsilon^2}j. 
\end{split}
\right.
\end{equation}
where 
\begin{equation}
\label{r_integral}
\rho(t,x)=\int_0^1 r \,dv. 
\end{equation}
As $\varepsilon\rightarrow 0^+$, \eqref{RJ} yields
\[r=\rho, \qquad j=-\frac{v}{\sigma(x)}\partial_x \rho. \]
Substituting this into the first equation of (\ref{RJ}) and integrating over $v$, one gets the limiting diffusion equation \cite{BSS}: 
\begin{equation}
\label{Diff_LTE}
\left\{
\begin{split}
\displaystyle& j=-\frac{v}{\sigma(x)}\partial_x \rho, \\[4pt]
\displaystyle& \partial_t \rho =\partial_x \left(\frac{1}{3\sigma(x)} \partial_x\rho\right).
\end{split}
\right.
\end{equation}

We solve the diffusive relaxation system \eqref{RJ} by splitting it into a relaxation step, followed by a transport step. 
One can check details of the discretized scheme in \cite{JPT2}, we omit it here. 

%We remark that one only needs to solve $r$, $j$ for $v>0$ in the diffusive relaxation system \eqref{RJ}.
%To reconstruct $f$ for both positive and negative velocities, it is obvious that
%\begin{align*}
%&f(t,x,v) = r(t,x,v) + \varepsilon j(t,x,v),   \\[2pt]
%&f(t,x,-v) = r(t,x,v) - \varepsilon j(t,x,v), 
%\end{align*}
\hspace{2cm}

{\bf B. Weights Balance Algorithm }
We review the weight balance algorithm studied in \cite{wang2020understanding}, which designs appropriate weights of different loss terms in the total loss function.

\hspace{2cm}

\begin{algorithm}[H]
\SetAlgoLined
Consider a physics-informed neural network $f_{\theta}(x)$ with parameters $\theta$ and a loss function
$$
\mathcal{L} = \mathcal{L}_G+\sum_{i = 1}^M\lambda_i\mathcal{L}_i(\theta).
$$
where $\mathcal{L}_G$ is the governing equation loss and $\mathcal{L}_i$ are the other losses (initial condition and etc.). $\lambda_i$ are the weights to balance the interplay of the losses.\;
\For{$n = 1,...,S$}
{
Compute $\hat{\lambda}_i$ by:
$$
\hat{\lambda}_i = \frac{\max_{\theta}|\nabla_{\theta}\mathcal{L}_G(\theta_n)|}{\overline{|\nabla_{\theta}\mathcal{L}_i(\theta_n)|}},
$$
where $\overline{|\nabla_{\theta}\mathcal{L}_i(\theta_n)|}$ is the mean of $|\nabla_{\theta}\mathcal{L}_i(\theta_n)|$ with respected to $\theta_n$\;
Update the weights $\lambda_i$ using a moving average:
$$
\lambda_i = (1-\alpha)\lambda_i+\alpha\lambda_i, i = 1,...M,
$$
where $\alpha$ is a constant and the authors suggest that $\alpha = 0.9$\;
Update the parameter $\theta$ via the gradient descent\;
}
\caption{Learning rate annealing for the PINN \cite{wang2020understanding}}
\label{cite_algo}
\end{algorithm}

\bibliographystyle{abbrv}
\bibliography{references}

\end{document}